%% Plain TeX
%%

\newcount\secno
\newcount\prmno
% Les commandes
\def\section#1{\vskip1truecm
               \global\def\currenvir{section}
               \global\advance\secno by1\global\prmno=0
               {\bf \number\secno. {#1}}
               \smallskip}

\def\subsection{\global\def\currenvir{subsection}
                \global\advance\prmno by1
               \smallskip  \ind{ (\number\secno.\number\prmno) }}
\def\subsec{\global\def\currenvir{subsection}
                \global\advance\prmno by1\smallskip 
                { (\number\secno.\number\prmno)\ }}

\def\proclaim#1{\global\advance\prmno by 1
                {\bf #1 \the\secno.\the\prmno$.-$ }}

\long\def\th#1 \enonce#2\endth{%
   \medbreak\proclaim{#1}{\it #2}\global\def\currenvir{th}\smallskip}

\def\rem#1{\global\advance\prmno by 1
{\it #1} \the\secno.\the\prmno$.-$ }

\magnification 1250
\pretolerance=500 \tolerance=1000  \brokenpenalty=5000
\mathcode`A="7041 \mathcode`B="7042 \mathcode`C="7043
\mathcode`D="7044 \mathcode`E="7045 \mathcode`F="7046
\mathcode`G="7047 \mathcode`H="7048 \mathcode`I="7049
\mathcode`J="704A \mathcode`K="704B \mathcode`L="704C
\mathcode`M="704D \mathcode`N="704E \mathcode`O="704F
\mathcode`P="7050 \mathcode`Q="7051 \mathcode`R="7052
\mathcode`S="7053 \mathcode`T="7054 \mathcode`U="7055
\mathcode`V="7056 \mathcode`W="7057 \mathcode`X="7058
\mathcode`Y="7059 \mathcode`Z="705A
\def\spacedmath#1{\def\packedmath##1${\bgroup\mathsurround =0pt##1\egroup$}
\mathsurround#1
\everymath={\packedmath}\everydisplay={\mathsurround=0pt}}
 \spacedmath{2pt}
\def\qfl#1{\buildrel {#1}\over {\longrightarrow}}
\def\mono{\lhook\joinrel\mathrel{\longrightarrow}}
\def\iso{\vbox{\hbox to .8cm{\hfill{$\scriptstyle\sim$}\hfill}
\nointerlineskip\hbox to .8cm{{\hfill$\longrightarrow $\hfill}} }}
\def\sdir_#1^#2{\mathrel{\mathop{\kern0pt\oplus}\limits_{#1}^{#2}}}

\font\eightrm=cmr8
\font\sixrm=cmr6
\font\san=cmssdc10
\def\ext{\hbox{\san \char3}}
\def\sym{\hbox{\san \char83}}
\def\pc#1{\tenrm#1\sevenrm}
\def\tx{\kern-1.5pt -}
\def\cqfd{\kern 2truemm\unskip\penalty 500\vrule height 4pt depth 0pt width
4pt\medbreak} 
\def\no{n\up{o}\kern 2pt}
\def\ind{\par\hskip 1truecm\relax}
\def\indp{\par\hskip 0.5truecm\relax}
\def\moins{\mathrel{\hbox{\vrule height 3pt depth -2pt width 6pt}}}
\def\Hom{\mathop{\rm Hom}\nolimits}
\def\Ext{\mathop{\rm Ext}\nolimits}
\def\Aut{\mathop{\rm Aut}\nolimits}
\def\det{\mathop{\rm det}\nolimits}

\def\rk{\mathop{\rm rk\,}\nolimits}
\def\mult{\mathop{\rm mult}\nolimits}
\def\pot#1{^{\scriptscriptstyle\otimes#1 }}
\def\pop#1{^{\scriptscriptstyle\oplus#1 }}
\def\pv{ {\it Proof} :\ }
\def\oc{{\cal O}_C}

\def\P{\Bbb P}
\def\C{\Bbb C}

\def\G{\Bbb G}
\font\pal=cmsy7
\def\su#1{{\cal SU}_C(#1)}
\def\sp#1{{\cal S}\kern-1pt\raise-1pt\hbox{\pal P}^{}_C(#1)}

\frenchspacing
\input xy
\xyoption{all}
\input amssym.def
\input amssym
\vsize = 25truecm
\hsize = 16.1truecm
\voffset = -.5truecm
\parindent=0cm
\baselineskip15pt
\overfullrule=0pt

\centerline{\bf Vector bundles and theta functions
on  curves of genus 2 and 3}\smallskip 
\centerline{Arnaud {\pc BEAUVILLE}}
\vskip1truecm
{\bf Introduction}
\ind Let $C$ be a smooth projective curve, of genus $g\ge 2$.
The moduli space $\su{r}$ of semi-stable vector bundles of rank
$r$ on $C$, with trivial determinant, is a normal projective
variety, wich can be considered as a non-abelian analogue of the
Jacobian variety  $JC$. It is actually related to $JC$ by the
following construction, which goes back (at least) to [N-R].
\def\jg{J^{g-1}}
Let $\jg$ be the translate of $JC$ parameterizing line
bundles of degree $g-1$ on $C$, and $\Theta \i \jg$ the
canonical theta divisor. For $E\in \su{r}$, consider the
locus
$$\Theta _E:=\{L\in \jg\ |\ H^0(C,E\otimes L)\not= 0\}\
.$$Then either $\Theta _E=\jg$, or $\Theta _E$ is in a natural way a
divisor in
$\jg$, belonging to the linear system $|r\Theta |$. In this way
we get a rational map
$$\theta:\su{r}\dasharrow |r\Theta |$$which is the
most obvious rational map of $\su{r}$ in a projective space:
 it can be identified to the map $\varphi_{\cal L}:\su{r}\dasharrow
{\P}(H^0(\su{r},{\cal L})^*)$ given by the global sections of the determinant
bundle
${\cal L}$, the positive generator of the Picard group of $\su{r}$
[B-N-R].
\ind For $r=2$ the map  $\theta$ is  an 
embedding if  $C$ is not hyperelliptic  [vG-I]. We consider in this paper the
higher rank case, where very little is known. The first part is devoted to  the
case $g=2$. There a curious numerical coincidence occurs, namely
$$\dim \su{r}=\dim |r\Theta |=r^2-1\ .$$
\ind For $r=2$ $\theta$ is an isomorphism [N-R]; for $r=3$ it is
a double covering, ramified along a sextic hypersurface which is the dual of
the ``Coble cubic" [O]. We will prove:
\smallskip 
{\bf Theorem A}$.-$ {\it For a curve $C$ of genus $2$, the map 
$\theta:\su{r}\dasharrow |r\Theta |$ is generically
finite} (or, equivalently, dominant). {\it It admits some  fibers of
dimension $\ge [{r\over 2}]-1$}.
\ind Our method is to consider the fibre of $\theta $ over a reducible element
of $|r\Theta |$ of the form $\Theta +\Delta $, where $\Delta $ is general in
$|(r-1)\Theta |$. The main point is to show that this fibre restricted to the
stable locus of $\su{r}$ is finite. The other elements of the fibre are the classes
of  the  bundles $\oc\oplus F$, with $\Theta _F=\Delta $; reasoning by
induction on $r$ we may assume that there are finitely many such $F$, and 
this gives the first assertion of the theorem (\S 1). The second one follows from
considering the restriction of $\theta $ to a particular class of vector bundles,
namely the symplectic bundles (\S 2).
\ind The method is not, in principle, restricted  to genus 2 curves -- but the
geometry in higher genus becomes much more intricate. In the second part
of the paper (\S 3) we will apply it to  rank 3 bundles in genus 3. Our
result is:
\smallskip 
{\bf Theorem B}$.-$ {\it Let $C$ be a curve of genus $3$. The map $\theta
:\su{3}\rightarrow |3\Theta |$ is a finite morphism}.
\ind This means that a semi-stable vector bundle of rank 3 on $C$ has always a
theta divisor; or alternatively
(see e.g. [B1]), that the linear system $|{\cal L}|$ on $\su{3}$ is base point free.
\ind This is not a big surprise since the result is already known for a {\it
generic} curve of genus 3 [R]. We believe, however, that the method is more
interesting than the result itself. In fact we translate the problem into an
elementary  question of projective geometry: what are the continuous families
of planes in $\P^5$ such that any two planes of the family intersect? It turns
out that this question has been completely (and beautifully) solved by Morin
[M]. Translating back his result into the language of vector bundles we get  a
complete list of the stable rank 3 bundles $E$ of degree 0 such that $\Theta
_E\supset \Theta $ (Theorem
3.1 below). Theorem B follows as a corollary. 
{\bigskip {\eightrm\baselineskip=12pt
\leftskip1cm\rightskip1cm\hskip0.8truecm  I am very much indebted to  C.
Ciliberto for pointing out the paper of Morin and for making it
accessible to me.\par}}\bigskip
{\it Notations} : \ind Throughout the paper we will work with a complex curve 
$C$ (smooth, projective, connected), of genus $g$. If
$E$ is a vector bundle on $C$, we will write $H^0(E)$ for 
$H^0(C,E)$, and $h^0(E)$ for its dimension.

\section{Genus 2: the generic finiteness}
\ind In this section we assume $g=2$. The first part of  theorem A follows from
a slightly more precise result:
\th Proposition
\enonce Let $\Delta $ be a general divisor in $|(r-1)\Theta |$.
The fibre $\theta ^{-1} (\Theta +\Delta )$ is finite and non-empty.
\endth
\subsection We will prove the Proposition by induction on $r$. Let
$[E_0]\in\theta ^{-1} (\Theta +\Delta )$. If it is
 not stable, it is the class of a direct sum
$\sdir_{i}^{}E_i$, so that   $\Theta _{E_0}=\sum_i\Theta _{E_i}$; thus
$[E_0]$ is the class in $\su{r}$ of $\oc\oplus F$ for some $F\in \su{r-1}$ with
$\Theta _F=\Delta
$. By the induction hypothesis there exists only finitely many such $F$, and
there exists at least one.
\ind Thus we can  assume that $E_0$ is stable. Let $E:=E^*_0\otimes
K_C$. We have $h^0(E)=r$ by Riemann-Roch
and the stability of $E_0$. The inclusion $\Theta =C\i \Theta _{E_0}$
means that $h^0(E_0(p))\ge 1$ for all $p\in C$, or equivalently by Serre
duality $h^0(E(-p))\ge 1$; this implies that  the subsheaf $F$ of
$E$  generated by the global sections of $E$ is of  rank 
$<r$.
Moreover if $p$ does not belong to 
$\Delta $, it is a smooth point of $\Theta _{E_0}$, and thus satisfies 
$h^0(E(-p))= 1$ (see e.g. [L], \S V); therefore $\rk F=r-1$ (otherwise
we would have $h^0(E(-p))\ge h^0(F(-p))\ge 2$).
\subsection   Let $Z$ be a component of the locus of stable bundles
$E$ of rank $r$ and determinant $K\pot{r}$, with the  property that $H^0(E)$
span a sub-bundle of rank $r-1$ of $E$. We will prove the inequality
$\dim Z\le \dim |(r-1)\Theta |$. It implies  that the general fibre of
$\theta:Z\dasharrow \Theta +|(r-1)\Theta |$ is
finite (possibly empty), so the Proposition follows.
\ind Let $E$ be a general element of $Z$, and let
 $F$ be the sub-bundle of $E$ spanned by $H^0(E)$. 
Put
$L:=\det F$ and $d=\deg 
F=\deg L$; we have an  exact sequence
$$0\rightarrow L^{-1} \longrightarrow H^0(E)\otimes _{\C}{\cal
O}_C\longrightarrow F\rightarrow 0\ ,$$hence a linear map
$H^0(E)^*\rightarrow H^0(L)$. 
 Let $s=r-\dim H^0(C,F^*)$ be the rank of that map. Then $F=\oc^{r-s}\oplus
G$, where $G$ is a vector bundle of rank $s-1$ with $h^0(G)=s$,
$h^0(G^*)=0$, which fits into  an exact sequence
$$0\rightarrow L^{-1} \longrightarrow {\cal
O}_C^{s}\longrightarrow G\rightarrow 0\ .$$
\ind The quotient ${\cal M}=E/F$ is 
the direct sum of a line bundle $M$
and a torsion sheaf ${\cal T}$. We have
$c_1(M)+c_1({\cal T})=rc_1(K_C)-c_1(L)$, and this formula
determines
$M$ once ${\cal T}$ and $L$ are given. We denote by $t$ the
length of 
${\cal T}$. 

\subsection To summarize, we have associated to a general bundle
$E$ in
$Z$ integers $s,d,t$ and
\ind -- a line bundle $L$ of degree $d$, and a $s$\tx dimensional
subspace  $V\i H^0(C,L)$ generating $L$; from these data we
define
$G$ as the cokernel of the natural map $L^{-1} \rightarrow
V^*\otimes \oc$, and put
$F:=\oc^{r-s}\oplus G$;
\ind -- a torsion sheaf ${\cal T}$ of length $t$ and an extension
$$0\rightarrow F \longrightarrow E\longrightarrow
 M\oplus {\cal T}\rightarrow 0\ ,\eqno({\cal E})$$
where the line bundle $M$ is determined by
$c_1(M)=rc_1(K_C)-c_1(L)-c_1({\cal T})$.
\ind The integers $s,d,t$ are bounded: we have $s\le
r$, $t\le 2r-d$, and $d<2(r-1)$ by the stability of $E$. Observe
also that $d\ge 3$: indeed $L$ is generated by its global sections,
and cannot be isomorphic to   $K_C$ since otherwise
$F$ would contain a copy of $K_C$,  contradicting the stability
of $E$.
\ind The data $(L,V,{\cal T},{\cal E})$ are parameterized by a
variety dominating 
$Z$; we will bound its dimension. The line bundle $L$ depends on
2  parameters.  We have $h^0(L)=d-1$ since $d\ge 3$, therefore
the subspace $V\i H^0(L)$ depends on $s(d-1-s)$ parameters. The
torsion sheaf ${\cal T}$ depends on $t$ parameters. Over the
variety parameterizing these data we build a vector bundle with
fibre $\Ext^1({\cal M}, F)$, with ${\cal M}=M\oplus{\cal T}$, $M$
and $F$ being determined as above. The group
$\Aut({\cal M})\times \Aut(F)$ acts on $\Ext^1({\cal M},
F)$, with the group $\C^*$ of homotheties of ${\cal M}$ and
$F$ acting in the same way; in fact, since the middle term of the
extensions we are interested in is stable, the stabilizer of  a general
extension class is
$\C^*$. This gives a bound
 $$\dim Z \le  2+s(d-1-s)+t+\dim\Ext^1({\cal
M},F)-\dim\Aut({\cal M})-\dim\Aut(F)+1\ .$$
\ind Let us estimate the  dimensions which appear in the right
hand side. We have $\Hom(M,F)\i \Hom(M,E)=0$
because
$E$ is stable, hence by Riemann-Roch
$$\dim \Ext^1(M\oplus {\cal T}, F)=(r-1)(2r+1)-dr\ .$$ 
\ind  The group
$\Aut(F)=\Aut(\oc^{r-s}\oplus G)$ contains the group of matrices
$\pmatrix{u & 0\cr v&\lambda}$, with $u\in \Aut(\oc^{r-s})$,
$v\in
\Hom(\oc^{r-s}, G)$, $\lambda\in\C^*$; this group has
dimension $$(r-s)^2+s(r-s)+1=r(r-s)+1\ .$$ The group $\Aut({\cal
T})$ has dimension at least $t$, so similarly $\Aut({\cal M})$ has
dimension $\ge 2t+1$. We get finally:
$$\eqalign{\dim Z&\le 2+s(d-1-s)+t+(r-1)(2r+1)-dr-r(r-s)-2t-1
\cr
&= (r-1)^2-1-(d-1-s)(r-s)-t}$$
Since $d-1=h^0(L)\ge s$, this implies
$\dim Z\le (r-1)^2 -1=\dim|(r-1)\Theta |$ as required.\cqfd

\section{Symplectic bundles}
\ind Let $C$ be a curve of genus $g\ge 2$, and $r$ a  positive
integer. The moduli space $\su{r}$ has a natural involution
$D:E\mapsto E^*$. Let $\iota$ be the involution $L\mapsto
K_C\otimes L^{-1}$ of $J^{g-1}$. The diagram  $$\xymatrix 
@M=6pt{\su{r} \ar[r]^{D }\ar@{-->}[d]_{\theta }& 
\su{r}\ar@{-->}[d]^{\theta}\\ |r\Theta |\ar[r]^{\iota^*}& |r\Theta|}$$
is commutative.
\ind Assume  now  that $r$ is even. Let $\sp{r}$  be the moduli
space of semi-stable symplectic  bundles of rank $r$ on
$C$. This is  a normal connected projective variety, with a
forgetful morphism to $\su{r}$, which is an embedding on the stable
locus. It is contained in the fixed locus of $D$, thus its image
under $\theta $ is contained in the fixed locus of $\iota ^*$.  
\ind This fixed locus is described for instance in [B-L], ch. 4, \S 6
(up to a translation from $JC$ to $J^{g-1}$). The involution $\iota^*$
acts linearly on $|r\Theta |$ and has 2 fixed spaces $|r\Theta |^+$
and $|r\Theta |^-$:  a symmetric divisor  in
$|r\Theta |$ is in $|r\Theta |^+$ (resp. $|r\Theta |^-$) if and only if
its multiplicity at any theta-characteristic $\kappa \in J^{g-1}$ is
even (resp. odd).  
The dimension of $|r\Theta |^{\pm}$ is ${1\over
2}(r^g\pm 2^g)-1$. 

\th Proposition
\enonce  $\theta:\su{r}\dasharrow|r\Theta |$  induces a rational
map from $\sp{r}$  to $|r\Theta |^+$. \endth
\pv Since $\sp{r}$ is connected, it
suffices to find one semi-stable bundle $E$ which admits a
symplectic form, and such that $\Theta _E\in |r\Theta|^+$.  We take
$E=F\oplus F^*$ with the standard alternate form, where $F\in\su{r/2}$ 
admits a theta divisor. Then
$\Theta_E=\Theta_F+\iota^*\Theta_F$. Thus if  $\kappa  \in
J^{g-1}$ is a  theta-characteristic, we have $\mult_\kappa (\Theta
_E)=2\mult_\kappa (\Theta _F)$, hence  $\Theta _E\in
|r\Theta|^+$.\cqfd   \ind Let us go back to the case $g=2$.
\th Proposition \enonce If $C$ has genus $2$, some fibres of $\theta
:\su{r}\dasharrow   |r\Theta|$ have dimension $\ge [{r\over 2}]-1$.
\endth
\pv If $r$ is even, $\theta $ induces  a rational map $\theta_{sp}
:\sp{r}\dasharrow  |r\Theta|^+$ (Prop. 2.1). We have 
$$\dim \sp{r}={1\over 2}r(r+1)\quad ,\quad \dim  |r\Theta|^+={r^2\over
2}+1\ ,$$hence  the fibres have dimension $\ge {r\over 2}-1$.
\ind If $r$ is odd, consider the bundle  $E\oplus\oc$, for $E$
general in $\sp{r-1}$; by what we have just seen $\theta $ is
defined at
$E$, and its fibre at $E$  has dimension $\ge {r-1\over 2}-1$.\cqfd 
\rem{Remark}The {\it degree} of
$\theta_r:\su{r}\dasharrow |r\Theta |$ grows exponentially with $r$: indeed
the commutative diagram
$$\xymatrix 
@M=6pt{\su{r}\times \su{s} \ar[r]^>>>>>{\oplus
}\ar@{-->}[d]_{\theta_r\times
\theta_s }& 
\su{r+s}\ar@{-->}[d]^{\theta_{r+s}}\\ 
|r\Theta |\times |s\Theta |\ar[r]^{+}& |(r+s)\Theta|}$$
shows that $\deg \theta_{r+s}\ge \deg \theta_r \cdot \deg \theta_s$. Since
$\deg \theta_3=2$, we obtain $\deg \theta_r\ge 2^{[r/3]}$ (we expect the
actual value to be much higher). 
\section{Genus 3, rank 3}
\ind Recall that if $L$ is a line bundle on $C$ generated by its global sections,
the {\it evaluation bundle} $Q_L$ is defined through the exact sequence
$$0\rightarrow Q_L^*\longrightarrow H^0(L)\otimes_{\C}\oc\longrightarrow
L\rightarrow 0\ ;$$it has rank $h^0(L)-1$ and determinant $L$. 
\th Theorem
\enonce Let $C$ be a curve of genus $3$, and $E_0$  a stable vector bundle
of rank $3$ and degree $0$ on $C$, such that $\Theta _{E_0}\supset\Theta $.
Then $C$ is not hyperelliptic, and $E_0$ is one of the following bundles:
\indp {\rm a)} The vector bundles $E_N:=Q_{K\otimes N}\otimes N^{-1} $,
for
$N\in J^2\moins\Theta
$;
\indp{\rm b)}  The vector bundle ${\cal E}nd_0(Q_K)$ of traceless
endomorphisms of $Q_K$.
\ind Conversely, the bundles in {\rm a)} and {\rm b)} are stable and admit a
theta divisor which contains $\Theta $.
\endth
\ind Thus all vector bundles in $\su{3}$ have a theta divisor; in other words, the
map $\theta :\su{3}\rightarrow |3\Theta |$ is a morphism. Since $\theta ^*{\cal
O}(1)={\cal L}$ is  ample, this morphism is finite: this implies Theorem B of the
introduction.\smallskip 
\ind The proof of Theorem 3.1 will occupy the rest of this section. Let
$E_0$ be a stable bundle in $\su{3}$ with $\Theta _{E_0}\supset\Theta $. We
will deal mainly with the bundle $E:=E^*\otimes K_C$.  It has  slope
4, degree 12 and satisfies
$h^1(E)=h^0(E_0)=0$ by stability of $E_0$, so that
$h^0(E)=6$ by Riemann-Roch. We first establish some properties of $E$ that will
be needed later on.
\th Lemma
\enonce Any rank $2$ sub-bundle $F$ of $E$ satisfies $h^0(F)\le 4$.
\endth
\pv Assume $h^0(F)\ge 5$. Let $A$ be a sub-line bundle of
$F$ of maximal degree; this degree is 
$\ge 2$ (since $h^0(F(-p-q)\ge 1$ for $p,q\in C$) and $\le 3$ by the
stability of $E$. Let
$B:=F/A$; again by stability of $E$ we have $\deg(F)\le 7$, hence 
 $\deg(B)\le 5$. Therefore  $$h^0(F)\le
h^0(A)+h^0(B)\le 2+3=5\ ;$$if equality holds, we have $h^0(A)=2$,
$h^0(B)=3$; moreover the  class of the extension
$$0\rightarrow A\longrightarrow F\longrightarrow B\rightarrow
0$$must be non-zero (because $E$ cannot contain a line bundle of
degree $\ge 4$), but must go to zero under the canonical map
$$\Ext^1(B,A)\longrightarrow \Hom(H^0(B),H^1(A))\ .$$In
particular this map cannot be injective; equivalently its transpose,  the
multiplication map 
$$ H^0(K\otimes A^{-1} )\otimes H^0(B)\longrightarrow
H^0(K\otimes A^{-1}\otimes B)$$cannot be surjective.
Now we distinguish two cases:
\ind a) If $\deg(A)=3$, we must have $A=K_C(-p)$ for some $p\in C$,
and $B=K_C$. But then the multiplication map  
$H^0(\oc(p) )\otimes H^0(K_C)\iso H^0(K_C(p))$ is an isomorphism.

\ind b) If $\deg(A)=2$,  $C$ is hyperelliptic
and $A$ is the hyperelliptic line bundle on $C$ (that is, $h^0(A)=\deg A=2$). If
$B=K_C$, the multiplication map 
$H^0(A)\otimes H^0(K_C)\rightarrow H^0(A\otimes K_C)$ is
surjective. So we must have $\deg(B)=5$. By the base point free pencil
trick, the multiplication map
$H^0(A)\otimes H^0(B)\rightarrow H^0(A\otimes B)$ is surjective if
and only if $H^1(B\otimes A^{-1} )=0$, that is, $H^0(K\otimes A\otimes
B^{-1} )=0$. This fails only if  $B\cong K(q)$ for some $q\in
C$. But in that case $B$, and therefore also $F$, are  not globally
generated. The subsheaf $F'$ of $F$ spanned by $H^0(F)$ has
$h^0(F')=5$, $\deg(F')\le 6$, and this is impossible by the
previous analysis.\cqfd
\th Lemma
\enonce Let $p,q$ be general points of $C$. Then $h^0(E(-p))=3$ and
$h^0(E(-p-q))=1$.
\endth
\pv If $h^0(E(-p))\ge 4$ for all $p\in C$, the global sections of $E$
span a sub-bundle $F$ of rank
$\le 2$ with  $h^0(F)=6$. This is impossible by Lemma
3.2. Similarly if
$h^0(E(-p-q))\ge 2$ for all $q$, the global sections of $E(-p)$ span a sub-line
bundle $L$ of
$E(-p)$ with $h^0(L)=3$, hence  $\deg L\ge 4$, contradicting
the stability of $E$.\cqfd
\ind Thus the spaces $\P(H^0(E(-p)))$ form a one-dimensional family of
planes in $\P(H^0(E))\cong \P^5$ with the property that any two of
them intersect. This situation has been thoroughly analyzed by Morin
[M]. \smallskip 

{\bf Theorem} (Morin)$.-$ {\it Any irreducible family of planes in
$\P^5$ such that any two planes of the family intersect is contained in
one of the following families:
\indp {\rm e1)} The planes passing trough a given point.
\indp {\rm e2)} The planes contained in a given hyperplane.
\indp {\rm e3)} The planes intersecting a given plane along a line.
\indp {\rm g1)} One of the family of generatrices of a smooth quadric in
$\P^5$.
\indp {\rm g2)} The family of planes cutting down a smooth conic on
the Veronese surface. 
\indp {\rm g3)} The family of planes in $\P^5$ tangent to the Veronese
surface}.
\subsec {\it The elementary cases}
\ind We will first show that our family of planes cannot satisfy one of
the elementary conditions e1) to e3).
\indp e1) This would mean that there exists a non-zero section
$s\in H^0(E)$ which vanishes at each point of $C$, a
contradiction.
\indp e2) In that case there exists a hyperplane $H$ in $H^0(E)$
such that $H^0(E(-p))\i H$ for all $p$ in $C$. It follows that $H$
span a sub-bundle $F$ of $E$ of rank $\le 2$, with $h^0(F)\ge 5$; this
contradicts Lemma 3.2.
\indp e3) In that case there exists a 3-dimensional subspace
$W$ in
$H^0(E)$ such that $\dim W\cap H^0(E(-p))\ge 2$ for all $p$ in
$C$. This implies that $W$ spans a sub-line bundle  $L$ of $E$
with $h^0(L)\ge 3$,  contradicting  the stability of $E$. \smallskip 
\subsec {\it The geometric cases}
\ind Suppose now that our family of planes $\P(H^0(E(-p)))\i \P(H^0(E))$ is
contained in one of the families  g1) to g3). We put $V:=H^0(E)$ and
consider the map $g:C\rightarrow \G(3,V)$ which associates to a
general  point $p$ of $C$ the subspace $H^0(E(-p))$ of $V$. This map
is defined by the sub-bundle $E'$ of $E$ spanned by $H^0(E)$; that is,
the universal exact sequence on $\G(3,V)$
$$0\rightarrow N \longrightarrow V\otimes {\cal
O}_{\G}\longrightarrow Q\rightarrow 0$$
pulls back to the exact sequence
$$0\rightarrow N_C \longrightarrow V\otimes \oc\longrightarrow
E'\rightarrow 0$$
on $C$,  where $(N_C)_p=H^0(E(-p))$ for $p$ general in $C$. The
Morin theorem tells us that $g$ factors as
$$g:C\qfl{f}\P^r\mono\G(3,V)\ ,$$where $r=2$ or 3 and  $\P^r$ is embedded
in $\G(3,V)$ as described in g1) to g3). Conversely if this holds, the vector
bundle $E'=g^*Q$ has the property that $h^0(E'(-p-q))\ge 1$ for all $p,q$ in
$C$.
\ind We will now analyze each of these cases and deduce from this the possibilities
for $E$. We put $L:=f^*{\cal O}_{\P^r}(1)$.
\smallskip 
 g1) {\it Planes in a quadric}
\ind Let $U$ be a 4-dimensional vector space, and $V=\ext^2U$. The equation
$v\wedge v=0$ for $v\in V$ defines a smooth quadric ${\cal Q}$ in $\P(V)$.  The
subvariety of $\G(3,V)$ parameterizing planes contained in ${\cal Q}$ has
two components, which are exchanged under the automorphism group of ${\cal
Q}$. One of these is the image of
  the map $\P^3=\P(U^*)\rightarrow \G(3,V)$ which maps the hyperplane $H\i
U$ to the  3-plane $\ext^2H\i\ext^2U$ $=V$. 
 The  Euler exact sequence 
$$0\rightarrow \Omega ^1_{\P^3}(1)\longrightarrow
U\otimes_{\C} {\cal O}_{\P^3} \longrightarrow  {\cal O}_{\P^3}
(1)\rightarrow 0$$ gives rise to an exact sequence
$$0\rightarrow \ext^2\bigl(\Omega ^1_{\P^3} (1)\bigr) \longrightarrow 
\ext^2U\otimes_{\C} {\cal O}_{\P^3} \longrightarrow
\Omega ^1_{\P^3}(2)\rightarrow 0$$
which is the pull back to $\P^3$ of the universal exact sequence on
$\G(3,V)$. 
\ind Thus $E'\cong f^*\Omega ^1_{\P^3}(2)$; the Euler exact sequence
twisted by ${\cal O}_{\P^3}(1)$ pulls back to
$$0\rightarrow E'\longrightarrow U\otimes_{\C}L\longrightarrow
L\pot{2}\rightarrow 0\ .$$
This implies $\det E'\cong L\pot{2}$,
hence $\deg L\le 6$. On the other hand the condition $h^0(E'(-p-q))\ge 1$ for
all $p,q$ in $C$ implies $h^0(L)\ge 3$ and therefore $\deg L\ge 4$. The map
$U\rightarrow H^0(L)$ must then be injective, because otherwise a copy of
$L$ would inject into $E'$, contradicting  the stability of $E$. This gives
$h^0(L)\ge 4$; the only possibility is
$\deg L=6$ and $h^0(L)=4$, hence $E'=E$ and
$U=H^0(L)$. Thus $E$  is isomorphic to $Q_L^*\otimes L$, where $Q_L$ is the
 evaluation bundle of $L$. This vector bundle  is analyzed in [B2]: it
always admits a theta divisor, and it is stable if and only if $C$ is not
hyperelliptic and 
$L$ is very ample, that is, $L=K_C\otimes N$ with $\deg N=2$, $h^0(N)=0$.
Dualizing we find $E_0=E_N:=Q_{K\otimes N}\otimes N^{-1}$; this gives 
case a) of the theorem.
\medskip 
 g2) {\it Secant planes to the Veronese surface}
\ind Let $U$ be a 3-dimensional vector space, and $V=\sym^2U$. The
Veronese surface $S$ is the image of the  map $u\mapsto u^2$ from
$\P(U)$ into $\P(V)$.  The family of planes which cut $S$ along a
conic is the image of the map
$\P^2=\P(U^*)\rightarrow \G(3,V)$ which maps a 2-plane 
$H\i U$ to
$\sym^2H\i \sym^2U$. The  pull back to $\P^2$ of the universal exact
sequence on $\G(3,V)$ is the sequence
 $$0\rightarrow \sym^2(\Omega ^1_{\P^2}(1))\longrightarrow
V\otimes _{\C}{\cal O}_{\P^2}\longrightarrow
{\cal O}_{\P^2}(1)^3\rightarrow 0$$obtained by taking the symmetric
square of the Euler exact sequence on $\P^2$.
\ind  Thus  $E'$ is
isomorphic to $L^{\oplus 3}$. Since $E$ is stable this
implies $\deg L\le 3$, while the inequality 
$h^0(E'(-p-q))\ge 1$ imposes
$h^0(L)\ge 3$, a contradiction. 
\medskip
 g3) {\it Tangent planes to the Veronese surface}
\ind  Consider again the Veronese surface $S$, image of the square
map $\P(U)\rightarrow \P(V)$.  The
projective tangent bundle of
$S$ in $\P(V)$ is $\P_S(\widetilde{T}_S)$, where $\widetilde{T}_S$ appears in
the extension
$$0\rightarrow {\cal O}_S\longrightarrow  \widetilde{T}_S\longrightarrow
T_S\rightarrow 0$$
with class $c_1({\cal O}_{\P(V)}(1)^{}_{|S})\in H^1(S,\Omega
^1_S)$; the Euler exact sequence provides an isomorphism
$\widetilde{T}_S\cong U\otimes {\cal O}_{\P(U)}(1)$. Similarly we have an extension
$\widetilde{T}_{\P(V)}$ of 
$T_{\P(V)}$ by ${\cal O}_{\P(V)}$ and an isomorphism
$\widetilde{T}_{\P(V)}\cong V\otimes {\cal O}_{\P(V)}(1)$. These bundles fit
into a normal exact sequence
$$0\rightarrow \widetilde{T}_S\longrightarrow
\widetilde{T}_{\P(V)}{}^{}_{|S}\longrightarrow N_{S/\P(V)}\rightarrow 0\
,$$that is, after a twist by ${\cal O}_S(-2)$,
$$0\rightarrow U\otimes
{\cal O}_{S}(-1)\longrightarrow
V\otimes {\cal O}_{S}\longrightarrow N_{S/\P(V)}(-2)\rightarrow 0\ ,$$
which is the pull back  to $S$ of the universal exact
sequence on $\G(3,V)$. Recall that the second
fundamental form gives an isomorphism $N_{S/\P^5}\cong
\sym^2T_{S}$ (see for instance [G-H]).  
\ind  Thus $E'=\sym^2f^*(T_{\P^2}(-1))$. This gives 
 $\det E'=L\pot{3}$, hence $\deg L\le
4$. On the other hand we have $h^0(L)\ge 3$:  otherwise the image of
$C$ in $\P^5$ is a conic  $c\i S$, and all tangent planes to $S$ along
$c$  meet the plane of $c$ along a line, so that we are in case e3). Therefore
 $L=K_C$, $E=E'$. The Euler
exact sequence shows that $f^*(T_{\P^2}(-1))$ is isomorphic to the
evaluation bundle $Q_K$ of $K_C$, so that $E\cong \sym^2Q_K$.  Using the
canonical isomorphism $\sym^2F\otimes (\det F)^{-1}\iso {\cal E}nd_0(F) $ 
for a rank 2 bundle $F$ we get $E_0\cong {\cal E}nd_0(Q_K)$.
\ind The vector bundles  $Q_K$, and therefore
${\cal E}nd_0(Q_K)$, are semi-stable.  If $C$ is  hyperelliptic, $Q_K$ is
isomorphic to $H\oplus H$, where $H$ is the hyperelliptic line bundle, hence
${\cal E}nd_0(Q_K)\cong \oc\pop{3}$. 
\ind Assume now that $C$ is not hyperelliptic; then  $Q_K$ is stable [P-R]. If 
${\cal E}:={\cal E}nd_0(Q_K)$ is not stable, it admits as sub- or quotient
bundle a line bundle  of degree 0; this means that there exists a non-zero
homomorphism $Q_K\rightarrow Q_K\otimes M$, with $M\in JC$, which 
 must be an isomorphism because $Q_K$ is stable. Taking determinants gives
$M\pot{2}\cong\oc$. Since $C$ is not hyperelliptic $M$ cannot be written
$\oc(p-q)$ with $p,q\in C$; therefore $h^0(Q_K\otimes M)=0$ [P-R], so that
$Q_K\otimes M$ cannot be isomorphic to $Q_K$.
\ind It remains to prove that ${\cal E}$ admits a theta divisor. What we have
proved so far is that ${\cal E}$ is the only stable rank 3 vector bundle of
degree 0 which might possibly satisfy $\Theta _{\cal E}=J^2$. But if this was
the case, all the vector bundles ${\cal E}\otimes M$, for $M\in JC$, should
have the same property -- an obvious contradiction.\cqfd
\smallskip 
\rem{Remarks}a) If we restrict ourselves to  $\su{3}$, we
find 37 stable bundles, namely ${\cal E}nd_0(Q_K)$ and the bundles
$E_\kappa$ where $\kappa$ is an even
theta-characteristic. These bundles appear already in [P], in a somewhat
disguised form: one can show indeed that $E_\kappa$ is isomorphic to ${\cal
E}nd_0(A(\kappa,L,x))$, where
$A(\kappa,L,x)$ is the {\it Aronhold bundle} defined in [P] (up to a twist, this
bundle depends only on $\kappa$).
\ind b) The theta divisor of $E_N$ is determined in
[B2]: it is equal to $\Theta +\Delta_N$, where $\Delta_N$ is the translate by
$N$ of the divisor $C-C$ in $JC$. The theta divisor of ${\cal E}nd_0(Q_K)$ is
$\Theta +\Xi$, where $\Xi$  is an interesting canonical
element of $|2\Theta |$. One can show that the trace of $\Xi$ on $\Theta
\cong \sym^2C$ is the locus of divisors $p+q$ such that the residual
intersection points of $C$ with the line $\langle p,q\rangle$ are harmonically
conjugate with respect to $p,q$ (here we view $C$ as a plane quartic). 
\ind c) Let $X\i |3\Theta |$ be the closed subvariety of divisors of the form 
$\Theta +\Theta _E$ for some $E$ in $\su{2}$. It follows from Theorem
3.1 and the above remarks that the fibre of $\theta :\su{3}\dasharrow
|3\Theta |$ over a general point of $X$ is reduced to one element, while $\theta
^{-1} (\Theta +\Delta _\kappa )$, for $\kappa $ an even theta-characteristic,
has 2 elements, namely $E_\kappa$ and $\oc\oplus (Q_K\otimes \kappa
^{-1})$. From general principles this implies that the variety $\theta (\su{3})$
is {\it  not normal} at the 36 points $\Theta +\Delta _\kappa$ (see for instance
[EGA], 15.5.3). 
\ind d) Assume that the N\'eron-Severi group of $JC$ has rank 1 -- this holds
if $C$ is general enough. Then a reducible divisor in $|3\Theta |$ must
contain a translate of $\Theta $. We thus deduce from Theorem 3.1 that
the stable vector bundles of rank 3 and degree 0 on $C$ which admit a
reducible theta divisor are those of the form $E_N\otimes M$ or ${\cal
E}nd_0(Q_K)\otimes M$, for $M\in JC$. 

\vskip2cm
\centerline{ REFERENCES} \vglue15pt\baselineskip12.8pt
\def\num#1{\smallskip\item{\hbox to\parindent{\enskip [#1]\hfill}}}
\parindent=1.38cm 
\num{B1} A. {\pc BEAUVILLE}: {\sl Vector bundles on curves and generalized
theta functions: recent results and open problems}. Current Topics in Complex
Algebraic Geometry, MSRI Publications {\bf 28}, 17-33; Cambridge University
Press (1995).
\num{B2} A. {\pc BEAUVILLE}: {\sl Some stable vector bundles  with reducible
theta divisors}. Manuscripta Math. {\bf 110}, 343--349 (2003).

\num{B-N-R} A. {\pc BEAUVILLE}, M.S. {\pc NARASIMHAN}, S. {\pc
RAMANAN}: {\sl Spectral curves and the generalised theta
divisor}. J. Reine Angew. Math. {\bf 398} (1989), 169--179. 

\num{B-L} C. {\pc BIRKENHAKE}, H. {\pc LANGE}: {\sl Complex abelian
varieties}. Grund.  Math. Wiss. {\bf 302}.
Springer-Verlag, Berlin, 1992. 
\num{EGA} A. {\pc GROTHENDIECK}: {\sl \'El\'ements de g\'eom\'etrie 
alg\'ebrique}, IV. Inst. Hautes \'Etudes Sci. Publ. Math. {\bf  32} (1967).
\num{vG-I} {\pc B. VAN} {\pc GEEMEN}, E. {\pc IZADI}: {\sl The tangent
space to the moduli space of vector bundles on a curve and the singular
locus of the theta divisor of the Jacobian}. J. Algebraic Geom. {\bf 10}
(2001),  133--177. 
\num{G-H} P. {\pc GRIFFITHS}, J. {\pc HARRIS}: {\sl Algebraic geometry and
local differential geometry}. Ann. Sci. \'Ecole Norm. Sup. (4) {\bf 12} (1979),
 355--452.
\num{L} Y. {\pc LASZLO}: {\sl Un th\'eor\`eme de Riemann pour les
diviseurs th\^eta sur les espaces de modules de fibr\'es stables sur une
courbe}. Duke Math. J. {\bf 64} (1991),  333--347. 
\num{M} U. {\pc MORIN}: {\sl Sui sistemi di piani a due a due incidenti}. Atti
Ist. Veneto {\bf 89} (1930), 907--926.

\num{N-R} M.S. {\pc NARASIMHAN}, S. {\pc
RAMANAN}: {\sl Moduli of vector bundles on a compact Riemann
surface}. Ann. of Math. (2) {\bf 89} (1969), 14--51. 

\num{O} A. {\pc ORTEGA}:  {\sl On the moduli space of rank 3 vector
bundles on a genus 2 curve and the Coble cubic}. J. Algebraic Geom., to
appear.

\num{P-R} K. {\pc PARANJAPE}, S. {\pc
RAMANAN}: {\sl On the canonical ring of a curve}. 
Algebraic geometry and commutative algebra, Vol. II, 503--516, 
Kinokuniya, Tokyo (1988).
\num{P} C. {\pc PAULY}: {\sl Self-duality of Coble's quartic hypersurface and
applications}. Michigan Math. J. {\bf 50} (2002),  551--574.
\num{R} M. {\pc RAYNAUD}: {\sl Sections des fibr\'es vectoriels sur
une courbe}. Bull. Soc. Math. France {\bf 110} (1982),
103--125. 

\vskip1cm
\def\pc#1{\eightrm#1\sixrm}
\hfill\vtop{\eightrm\hbox to 5cm{\hfill Arnaud {\pc BEAUVILLE}\hfill}
 \hbox to 5cm{\hfill Institut Universitaire de France\hfill}\vskip-2pt
\hbox to 5cm{\hfill \&\hfill}\vskip-2pt
 \hbox to 5cm{\hfill Laboratoire J.-A. Dieudonn\'e\hfill}
 \hbox to 5cm{\sixrm\hfill UMR 6621 du CNRS\hfill}
\hbox to 5cm{\hfill {\pc UNIVERSIT\'E DE}  {\pc NICE}\hfill}
\hbox to 5cm{\hfill  Parc Valrose\hfill}
\hbox to 5cm{\hfill F-06108 {\pc NICE} Cedex 02\hfill}}
\end